\documentclass[a4paper, 11pt, twoside]{amsart}

\usepackage{amsmath}
\usepackage{amssymb}
\usepackage{amsfonts}
\usepackage[francais]{babel}

\def\R{{\mathbb R}}
\newcommand{\M}{{\mathcal{M}}}
\title[Pincement des sous-vari\'{e}t\'{e}s extrins\`{e}quement homog\`{e}nes]{Pincement des sous-vari\'{e}t\'{e}s extrins\`{e}quement homog\`{e}nes dans un espace euclidien}
\author[Peter Quast]{Peter Quast}
    \address{D\'{e}partement de math\'{e}matiques, Universit\'e de Fribourg, Suisse}
        \email{peter.quast@unifr.ch}
    \curraddr{Institut f\"{u}r Mathematik, Universit\"{a}t Augsburg, Allemagne}
    \email{peter.quast@math.uni-augsburg.de}
    \thanks{soutenu par le Fonds
national Suisse (FNS) sous le num\'{e}ro de projet PBFR--106367}

\date{22 d\'{e}cembre 2004}

\begin{document}

\maketitle

\begin{abstract}
\noindent Soit une vari\'{e}t\'{e} ferm\'{e}e $M$ immerg\'{e}e dans $\R^m.$
Supposons que le fibr\'{e} trivial $M\times \R^m=TM\oplus \nu M$ soit
muni d'une connexion presque m\'{e}trique $\tilde{\nabla}$ pr\'{e}servant
presque la d\'{e}composition de $M\times\R^m$ en fibr\'{e} tangent et
normal. Supposons de plus que, par rapport \`{a} cette connexion, la
diff\'{e}rence $\Gamma=\partial-\tilde{\nabla}$ \`{a} la d\'{e}riv\'{e}e usuelle
$\partial$ dans $\R^m$ soit presque parall\`{e}le. Alors $M$ admet une
immersion dans $\R^m$ en tant que sous-vari\'{e}t\'{e} extrins\`{e}quement
homog\`{e}ne.
\end{abstract}

\section*{Introduction} Dans le cas des vari\'{e}t\'{e}s riemanniennes
abstraites, une propri\'{e}t\'{e} g\'{e}om\'{e}trique est souvent localement
caract\'{e}ris\'{e}e par le parall\'{e}lisme d'un certain tenseur. Par exemple
\'{E}.~Cartan a caract\'{e}ris\'{e} les espace localement sym\'{e}triques par le
parall\'{e}lisme de la courbure riemannienne. Plus g\'{e}n\'{e}ralement, Nomizu
ainsi que Ambrose et Singer ont montr\'{e} qu'un espace localement
homog\`{e}ne poss\`{e}de une connexion m\'{e}trique $\tilde{\nabla}$ par rapport
\`{a} laquelle la courbure riemannienne et la diff\'{e}rence entre la
connexion riemannienne et $\tilde{\nabla}$ sont parall\`{e}les. La
stabilit\'{e} de ces caract\'{e}risations a \'{e}t\'{e} \'{e}tudi\'{e}e par Katsuda
(cf.~\cite{Kats-89}).

Si nous passons aux sous-vari\'{e}t\'{e}s compl\`{e}tes connexes de l'espace
euclidien, une technique due \`{a} Str\"{u}bing (cf.~\cite{Stru-79}) montre
que des propri\'{e}t\'{e}s globales de la g\'{e}om\'{e}trie extrins\`{e}que se laissent
souvent globalement caract\'{e}riser par le parall\'{e}lisme d'une certaine
structure. Une sous-vari\'{e}t\'{e} $M$ de $\R^m$ est appel\'{e}e
extrins\`{e}quement homog\`{e}ne, si pour toute paire de points $p$ et $q$
il existe une isom\'{e}trie de $\R^m$ qui envoie $p$ sur $q$ et qui
laisse $M$ invariant. Une sous-vari\'{e}t\'{e} extrins\`{e}quement homog\`{e}ne de
$\R^m$ est donc un orbite d'un sous-groupe du groupe des isom\'{e}tries
de $\R^m.$ Par analogie avec le th\'{e}or\`{e}me de Nomizu, Olmos
\cite{Olmo-93} et Eschenburg \cite{Esch-98} ont montr\'{e} la
caract\'{e}risation suivante des sous-vari\'{e}t\'{e}s extrins\`{e}quement homog\`{e}nes
de $\R^m:$ une sous-vari\'{e}t\'{e} ferm\'{e}e $M$ de $\R^m$ est extrins\`{e}quement
homog\`{e}ne, si et seulement si le fibr\'{e} trivial $M\times\R^m=TM\oplus
\nu M$ admet une connexion m\'{e}trique $\tilde{\nabla},$ appel\'{e}e {\it
connexion canonique,} satisfaisant
\begin{enumerate}
    \item $\tilde{\nabla}$ pr\'{e}serve la d\'{e}composition $M\times\R^m=TM\oplus
            \nu M,$
            c.\`{a}.d.~$\tilde{\nabla}=\tilde{\nabla}^T+\tilde{\nabla}^{\perp}.$
    \item La diff\'{e}rence $\Gamma=\partial-\tilde{\nabla}$ avec la d\'{e}riv\'{e}e usuelle $\partial$ de $\R^m$ est
            $\tilde{\nabla}$-parall\`{e}lle.
\end{enumerate}
Olmos et S\'{a}nchez \cite{Ol-Sa-91} ont montr\'{e} que si, de plus,
$\tilde{\nabla}^{\perp}$ co\"{\i}ncide avec la connexion normale
$\nabla^{\perp}$ induite par $\partial,$ et si $M$ est une
sous-vari\'{e}t\'{e} pleine, alors $M$ est une orbite d'une repr\'{e}sentation
d'isotropie d'un espace sym\'{e}trique semi-simple, aussi appel\'{e}e
repr\'{e}sentation de type $s,$ et vice-versa.\par

Si, de plus, la partie tangente de la connexion canonique est la
connexion riemannienne $\nabla,$ c.\`{a}.d. si
$\tilde{\nabla}=\nabla\oplus\nabla^{\perp},$ alors $M$ est
extrins\`{e}quement sym\'{e}trique (invariante par les r\'{e}flexions par
rapport aux espaces normaux) dans $\R^m$ (cf.~\cite{Feru-80}).

\section*{Le r\'{e}sultat} Soit $\M(\Lambda,d,n,m,\varepsilon)$
l'ensemble des triplets $(M^n,f,\tilde{\nabla})$ form\'{e}s d'une
vari\'{e}t\'{e} ferm\'{e}e $M$ de dimension $n,$ d'une immersion $f$ de $M$ dans
l'espace euclidien \`{a} $m$ dimensions et d'une connexion
$\tilde{\nabla}$ d\'{e}finie sur $M\times\R^m=TM\oplus \nu M,$ qui
satisfont aux conditions suivantes:
\begin{enumerate}
    \item Par rapport \`{a} la m\'{e}trique induite, le diam\`{e}tre de
    $M$ est major\'{e} par $d$ et la norme uniforme $||\alpha ||_0$ de la seconde forme fondamentale
    $\alpha$ est major\'{e}e par $\Lambda.$
    \item La connexion $\tilde{\nabla}$ poss\`{e}de les propri\'{e}t\'{e}s suivantes:
        \begin{enumerate}
            \item La norme uniforme du tenseur $\Gamma=\partial-\tilde{\nabla}$
                    est major\'{e}e par $\Lambda.$
            \item La norme uniforme des tenseurs suivants est major\'{e}e par $\varepsilon:$
                \begin{itemize}
                    \item $\tilde{\nabla} g,$ o\`{u} $g$ est la m\'{e}trique usuelle  sur $M\times\R^m,$
                    \item $\tilde{\alpha},$ o\`{u}
                    $\tilde{\alpha}(X,Y)=(\tilde{\nabla}_XY)^{\perp},\;
                    X,Y\in T_pM,\; p\in M,$
                    \item $\tilde{\nabla}\Gamma.$
                \end{itemize}
        \end{enumerate}
\end{enumerate}
Notons que les conditions mentionn\'{e}es au premier point donnent lieu
\`{a} une borne inf\'{e}rieure du rayon d'injectivit\'{e} et du volume de $M$
(cf.~\cite{Quas-04}) et excluent ainsi l'effondrement.\par

\paragraph*{Th\'{e}or\`{e}me (\cite{Quas-prep})} {\it Il existe une constante
$\varepsilon,$ d\'{e}pendant de
    $\Lambda,\,d,\,
    n$ et $m$ telle que, si un triplet $(M^n,f,\tilde{\nabla})$ est
    dans
    $\M(\Lambda,d,n,m,\varepsilon),$ alors $M$ admet une immersion
    dans $\R^m$ en tant que sous-vari\'{e}t\'{e} extrins\`{e}quement homog\`{e}ne.}

\subsection*{Aper\c{c}u de la d\'{e}monstration}

La d\'{e}monstration est inspir\'{e}e des m\'{e}thodes utilis\'{e}es par Katsuda
(cf.~\cite{Kats-89}) dans le cas des vari\'{e}t\'{e}s homog\`{e}nes abstraites.
Supposons par contraposition que pour chaque entier positif $i,$ il
existe un triplet $(M_i,f_i,\tilde{\nabla}_i)$ dans
$\M\left(\Lambda,d,n,m,{1\over i}\right),$ tel que les $M_i$ ne se
r\'{e}alisent pas comme sous-vari\'{e}t\'{e} extrins\`{e}quement homog\`{e}ne dans
$\R^m.$ \`{A} l'aide du th\'{e}or\`{e}me de convergence de Gromov et du th\'{e}or\`{e}me
d'inclusion des espaces de H\"{o}lder, nous construisons une vari\'{e}t\'{e}
limite $M$ et une immersion $f$ de classe $\mathcal{C}^1$ de $M$
dans $\R^m.$ Le th\'{e}or\`{e}me d'Arzel\`{a}-Ascoli nous donne une connexion
continue limite $\tilde{\nabla}$ sur $M\times \R^m.$ La perte de
r\'{e}gularit\'{e} lors du passage \`{a} la limite ne nous permet pas de
conclure directement que $\tilde{\nabla}$ est une connexion
canonique. Pour relever ce d\'{e}fi, nous construisons, \'{e}galement \`{a}
l'aide d'un processus limite, les isometries n\'{e}cessaires en
utilisant le transport parall\`{e}le. Ceci nous permet de d\'{e}duire que
$f(M)$ est extrins\`{e}quement homog\`{e}ne. Une d\'{e}monstration d\'{e}taill\'{e}e se
trouve dans \cite{Quas-prep}.

\subsubsection*{Remarque} Bien que $f$ ne soit \`{a} priori que
de classe $\mathcal{C}^1,\, f(M),$ \'{e}tant extrins\`{e}quement homog\`{e}ne,
est une sous-vari\'{e}t\'{e} lisse de $\R^m.$

\subsection*{Corollaires} Ce th\'{e}or\`{e}me admet des corollaires
analogues dans le cas des orbites des repr\'{e}\-sentations de type $s$
ainsi que dans le cas des sous-vari\'{e}t\'{e}s extrins\`{e}quement sym\'{e}triques.
Dans le premier cas, nous montrons que l'immersion limite $f$ est
m\^{e}me de classe $\mathcal{C}^2$ et donne donc lieu \`{a} une connexion
normale usuel $\nabla^{\perp},$ qui est continue. De plus, nous
utilisons le fait que les orbites pleines des repr\'{e}sentations de
type $s$ sont les vari\'{e}t\'{e}s dont le transport parall\`{e}le normal le
long d'une courbe est donn\'{e} par une isom\'{e}trie euclidienne qui laisse
$f(M)$ invariante (cf.~\cite{Ol-Sa-91}).

Le th\'{e}or\`{e}me analogue pour les sous-vari\'{e}t\'{e}s extrins\`{e}quement
sym\'{e}triques des espaces euclidiens se trouve dans \cite{Quas-04}.
Dans ce cas, on ne suppose plus l'existence d'une connexion
suppl\'{e}mentaire  $\tilde{\nabla},$ mais uniquement que la seconde
forme fondamentale soit presque parall\`{e}le. De plus, on ne doit plus
construire des isom\'{e}tries, car elles sont donn\'{e}es par les r\'{e}flexions
par rapport aux espaces normaux.

\subsection*{Remerciements} J'adresse mes remerciements \`{a}
J.-H.~Eschenburg et \`{a} P.~Ghanaat pour leur aide et leurs suggestions
lors de la pr\'{e}paration de ce r\'{e}sultat.


\end{document}